\documentclass{article}
\usepackage{graphicx,amsmath,amsthm,amsopn,amstext,amscd,amsfonts,amssymb}
\usepackage{natbib}

\def \tr {\mathop{\rm tr}\nolimits}
\def \re {\mathop{\rm Re}\nolimits}
\def \im {\mathop{\rm Im}\nolimits}
\def \Vol {\mathop{\rm Vol}\nolimits}
\def \etr {\mathop{\rm etr}\nolimits}
\def \diag {\mathop{\rm diag}\nolimits}

\renewenvironment{abstract}
                 {\vspace{6pt}
                  \begin{center}
                  \begin{minipage}{5in}
                  \centerline{\textbf{Abstract}}
                  \noindent\ignorespaces
                 }
                 {\end{minipage}\end{center}}
\newtheorem{theorem}{\textbf{Theorem}}[section]
\newtheorem{corollary}{\textbf{Corollary}}[section]
\newtheorem{lemma}{\textbf{Lemma}}[section]
\theoremstyle{definition}
\newtheorem{definition}{\textbf{Definition}}[section]

\title{\Large \textbf{Shape theory via affine transformation: Some generalisations}}
\author{
  \textbf{Jos\'e A. D\'{\i}az-Garc\'{\i}a} \thanks{Corresponding author\newline
    }\\
  {\normalsize Department of Statistics and Computation} \\
  {\normalsize Universidad Aut\'onoma Agraria Antonio Narro}\\
  {\normalsize 25350 Buenavista, Saltillo, Coahuila, M\'exico} \\
  {\normalsize E-mail: jadiaz@uaaan.mx} \\[2ex]
  \textbf{Francisco J. Caro-Lopera} \\
  {\normalsize Department of Basic Sciences} \\
  {\normalsize Universidad de Medell\'{\i}n} \\
  {\normalsize Carrera 87 No.30-65, of. 5-103}\\
  {\normalsize Medell\'{\i}n, Colombia}\\
  {\normalsize E-mail: fjcaro@udem.edu.co}\\
}
\date{}
\begin{document}
\maketitle

\begin{abstract}
This work sets the statistical affine shape theory in the context of real normed division
algebras. The general densities apply for every field: real, complex, quaternion, octonion, and
for any noncentral and non-isotropic elliptical distribution; then the separated published
works about real and complex shape distributions can be obtained as corollaries by a suitable
selection of the field parameter and univariate integrals involving the generator elliptical
function. As a particular case, the complex normal affine density is derived and applied in
brain magnetic resonance scans of normal and schizophrenic patients.
\end{abstract}

\section{Introduction}\label{sec1}

The literature of matrix-variate distributions (real, complex, quaternion, octonion) tells us
about a great effort for obtaining separately topics that were noticed recently (\citet{dg:09},
\citet{dggj:09}, \citet{dggj:10}) can be derived under a general approach.

In fact, many models and techniques are explored first in the real case, and then their
extensions to the complex case are proposed joint with all the necessary mathematical tools for
their development. From  the last 60 years we can citate  hundreds of examples of these
extensions, see \citet{h:55} and \citet{j:64}, \citet{m:82} and \citet{k:65}, \citet{da:80} and
\citet{rva:05a}, among many others examples.

In the statistical shape theory context, \citet{DM98} gives an important summary of diverse
techniques  in the real case. By other hand, \citet{mdm:06} studied some of the topics in
\citet{DM98} for the complex case.

There a number of techniques in shape theory, we focus in this paper in  the affine approach.
\citet{gm:93}(see also \citet{dgr:03}) proposed an alternative system shape coordinates termed
configuration or affine coordinates, randomly indexed by a matrix multivariate gaussian
distribution. Then, \citet{cdg:09} extended this  theory by replacing the normality assumption
with a matrix multivariate elliptical law. Those works were studied in the real field, so if we
follow the tradition of the literature, we can expect an extension to the complex case, for
example, by studying the new jacobians, integrals and computations.

Instead of this,  we propose an unified approach for the statistical theory of shape, by
studying the real, complex, quaternion and octonion cases in a simultaneous way. As we shall
see in section \ref{sec2}, these four cases are formally termed, real normed division algebras.
However as  usual, this type of generalisation have a price, in this case we need some concepts
and notation from the abstract algebra.

For the sake of completeness,  the case of the octonions is considered, but is important to
highlight that many of the results for the octonion field, only can be conjectured, because,
many theoretical problems about these numbers remain open, see \citet{dm:99}. In fact, the
relevance of the octonions  for understanding the real world is not clear at present, see
\citet{b:02}.

Section \ref{sec2} reviews some definitions and notation on real normed division algebras,
also, some concepts and integral properties of Jack polynomials and generalised hypergeometric
function are given; then in section \ref{sec3} a Jacobian with respect to Lebesgue measure for
real normed division algebras is obtained and the main results follows as a consequence;
section \ref{sec4} gives the affine shape distribution for several particular elliptical laws;
and finally, section \ref{sec5} shows an application from the literature of shape in complex
case.

\section{Preliminary results}\label{sec2}

A detailed discussion of real normed division algebras may be found in \citet{b:02} and
\citet{gr:87}, and of Jack polynomials and hypergeometric functions in \citet{S:97},
\citet{gr:87} and \citet{KE:06}. For convenience, we shall introduce some notations, although
in general we adhere to standard notations.

For us a \textbf{vector space} shall always be a finite-dimensional module over the field of
real numbers. An \textbf{algebra} $\mathfrak{F}$ shall be a vector space that is equipped with
a bilinear map $m: \mathfrak{F} \times \mathfrak{F} \rightarrow \mathfrak{F}$ termed
\emph{multiplication} and a nonzero element $1 \in \mathfrak{F}$ termed the \emph{unit} such
that $m(1,a) = m(a,1) = 1$. As usual, we abbreviate $m(a,b) = ab$ as $ab$. We do not assume
$\mathfrak{F}$ associative. Given an algebra, we shall freely think of real numbers as elements
of this algebra via the map $\omega \mapsto \omega 1$.

An algebra $\mathfrak{F}$ is a \textbf{division algebra} if given $a, b \in \mathfrak{F}$ with
$ab=0$, then either $a=0$ or $b=0$. Equivalently, $\mathfrak{F}$ is a division algebra if the
operation of left and right multiplications by any nonzero element are invertible. A
\textbf{normed division algebra} is an algebra $\mathfrak{F}$ that is also a normed vector
space with $||ab|| = ||a||||b||$. This implies that $\mathfrak{F}$ is a division algebra and
that $||1|| = 1$.

There are exactly four real finite-dimensional normed division algebras: real numbers, complex
numbers, quaternions and octonions, these being denoted generically as $\mathfrak{F}$, see
\citet{b:02}. All division algebras have a real dimension of $1, 2, 4$ or $8$, respectively,
whose dimension is denoted by $\beta$, see \citet[Theorems 1, 2 and 3]{b:02}.

Let $\mathcal{L}^{\beta}_{m,n}$ be the linear space of all $n \times m$ matrices of rank $m
\leq n$ over $\mathfrak{F}$ with $m$ distinct positive singular values, where $\mathfrak{F}$
denotes a \emph{real finite-dimensional normed division algebra}. Let $\mathfrak{F}^{n \times
m}$ be the set of all $n \times m$ matrices over $\mathfrak{F}$. The dimension of
$\mathfrak{F}^{n \times m}$ over $\Re$ is $\beta mn$. Let $\mathbf{A} \in \mathfrak{F}^{n
\times m}$, then $\mathbf{A}^{*} = \overline{\mathbf{A}}^{T}$ denotes the usual conjugate
transpose.

The set of matrices $\mathbf{H}_{1} \in \mathfrak{F}^{n \times m}$ such that
$\mathbf{H}_{1}^{*}\mathbf{H}_{1} = \mathbf{I}_{m}$ is a manifold denoted ${\mathcal
V}_{m,n}^{\beta}$, termed \emph{Stiefel manifold} ($\mathbf{H}_{1}$ is also known as
\emph{semi-orthogonal} ($\beta = 1$), \emph{semi-unitary} ($\beta = 2$), \emph{semi-symplectic}
($\beta = 4$) and \emph{semi-exceptional type} ($\beta = 8$) matrices, see \citet{dm:99}). The
dimension of $\mathcal{V}_{m,n}^{\beta}$ over $\Re$ is $[\beta mn - m(m-1)\beta/2 -m]$. In
particular, ${\mathcal V}_{m,m}^{\beta}$ with dimension over $\Re$, $[m(m+1)\beta/2 - m]$, is
the maximal compact subgroup $\mathfrak{U}^{\beta}(m)$ of ${\mathcal L}^{\beta}_{m,m}$ and
consist of all matrices $\mathbf{H} \in \mathfrak{F}^{m \times m}$ such that
$\mathbf{H}^{*}\mathbf{H} = \mathbf{I}_{m}$. Therefore, $\mathfrak{U}^{\beta}(m)$ is the
\emph{real orthogonal group} $\mathcal{O}(m)$ ($\beta = 1$), the \emph{unitary group}
$\mathcal{U}(m)$ ($\beta = 2$), \emph{compact symplectic group} $\mathcal{S}p(m)$ ($\beta = 4$)
or \emph{exceptional type matrices} $\mathcal{O}o(m)$ ($\beta = 8$), for $\mathfrak{F} = \Re$,
$\mathfrak{C}$, $\mathfrak{H}$ or $\mathfrak{O}$, respectively.

Denote by ${\mathfrak S}_{m}^{\beta}$ the real vector space of all $\mathbf{S} \in
\mathfrak{F}^{m \times m}$ such that $\mathbf{S} = \mathbf{S}^{*}$. Let
$\mathfrak{P}_{m}^{\beta}$ be the \emph{cone of positive definite matrices} $\mathbf{S} \in
\mathfrak{F}^{m \times m}$; then $\mathfrak{P}_{m}^{\beta}$ is an open subset of ${\mathfrak
S}_{m}^{\beta}$. Over $\Re$, ${\mathfrak S}_{m}^{\beta}$ consist of \emph{symmetric} matrices;
over $\mathfrak{C}$, \emph{Hermitian} matrices; over $\mathfrak{H}$, \emph{quaternionic
Hermitian} matrices (also termed \emph{self-dual matrices}) and over $\mathfrak{O}$,
\emph{octonionic Hermitian} matrices. Generically, the elements of $\mathfrak{S}_{m}^{\beta}$
are termed as \textbf{Hermitian matrices}, irrespective of the nature of $\mathfrak{F}$. The
dimension of $\mathfrak{S}_{m}^{\beta}$ over $\Re$ is $[m(m-1)\beta+2]/2$.

\noindent Let $\mathfrak{D}_{m}^{\beta}$ be the \emph{diagonal subgroup} of
$\mathcal{L}_{m,m}^{\beta}$ consisting of all $\mathbf{D} \in \mathfrak{F}^{m \times m}$,
$\mathbf{D} = \diag(d_{1}, \dots,d_{m})$.

For any matrix $\mathbf{X} \in \mathfrak{F}^{n \times m}$, $d\mathbf{X}$ denotes the\emph{
matrix of differentials} $(dx_{ij})$. Finally, we define the \emph{measure} or volume element
$(d\mathbf{X})$ when $\mathbf{X} \in \mathfrak{F}^{m \times n}, \mathfrak{S}_{m}^{\beta}$,
$\mathfrak{D}_{m}^{\beta}$ or $\mathcal{V}_{m,n}^{\beta}$, see \citet{d:02}.

If $\mathbf{X} \in \mathfrak{F}^{n \times m}$ then $(d\mathbf{X})$ (the Lebesgue measure in
$\mathfrak{F}^{n \times m}$) denotes the exterior product of the $\beta mn$ functionally
independent variables
$$
  (d\mathbf{X}) = \bigwedge_{i = 1}^{n}\bigwedge_{j = 1}^{m}dx_{ij} \quad \mbox{ where }
    \quad dx_{ij} = \bigwedge_{k = 1}^{\beta}dx_{ij}^{(k)}.
$$

If $\mathbf{S} \in \mathfrak{S}_{m}^{\beta}$ (or $\mathbf{S} \in \mathfrak{T}_{L}^{\beta}(m)$)
then $(d\mathbf{S})$ (the Lebesgue measure in $\mathfrak{S}_{m}^{\beta}$ or in
$\mathfrak{T}_{L}^{\beta}(m)$) denotes the exterior product of the $m(m+1)\beta/2$ functionally
independent variables (or denotes the exterior product of the $m(m-1)\beta/2 + n$ functionally
independent variables, if $s_{ii} \in \Re$ for all $i = 1, \dots, m$)
$$
  (d\mathbf{S}) = \left\{
                    \begin{array}{ll}
                      \displaystyle\bigwedge_{i \leq j}^{m}\bigwedge_{k = 1}^{\beta}ds_{ij}^{(k)}, &  \\
                      \displaystyle\bigwedge_{i=1}^{m} ds_{ii}\bigwedge_{i < j}^{m}\bigwedge_{k = 1}^{\beta}ds_{ij}^{(k)}, &
                       \hbox{if } s_{ii} \in \Re.
                    \end{array}
                  \right.
$$
Generally the context establishes the conditions on the elements of $\mathbf{S}$, that is, if
$s_{ij} \in \Re$, $\in \mathfrak{C}$, $\in \mathfrak{H}$ or $ \in \mathfrak{O}$. It is
considered that
$$
  (d\mathbf{S}) = \bigwedge_{i \leq j}^{m}\bigwedge_{k = 1}^{\beta}ds_{ij}^{(k)}
   \equiv \bigwedge_{i=1}^{m} ds_{ii}\bigwedge_{i < j}^{m}\bigwedge_{k =
1}^{\beta}ds_{ij}^{(k)}.
$$
Note that, the  Lebesgue measure $(d\mathbf{S})$ requires that $\mathbf{S} \in
\mathfrak{P}_{m}^{\beta}$, that is, $\mathbf{S}$ must be a non singular Hermitian matrix
(Hermitian definite positive matrix).

If $\mathbf{\Lambda} \in \mathfrak{D}_{m}^{\beta}$ then $(d\mathbf{\Lambda})$ (the Legesgue
measure in $\mathfrak{D}_{m}^{\beta}$) denotes the exterior product of the $\beta m$
functionally independent variables
$$
  (d\mathbf{\Lambda}) = \bigwedge_{i = 1}^{n}\bigwedge_{k = 1}^{\beta}d\lambda_{i}^{(k)}.
$$
If $\mathbf{H}_{1} \in \mathcal{V}_{m,n}^{\beta}$ then
$$
  (\mathbf{H}^{*}_{1}d\mathbf{H}_{1}) = \bigwedge_{i=1}^{n} \bigwedge_{j =i+1}^{m}
  \mathbf{h}_{j}^{*}d\mathbf{h}_{i}.
$$
where $\mathbf{H} = (\mathbf{H}_{1}|\mathbf{H}_{2}) = (\mathbf{h}_{1}, \dots,
\mathbf{h}_{m}|\mathbf{h}_{m+1}, \dots, \mathbf{h}_{n}) \in \mathfrak{U}^{\beta}(m)$. It can be
proved that this differential form does not depend on the choice of the matrix $\mathbf{H}_{2}$
When $m = 1$; $\mathcal{V}^{\beta}_{1,n}$ defines the unit sphere in $\mathfrak{F}^{n}$. This
is, of course, an $(n-1)\beta$- dimensional surface in $\mathfrak{F}^{n}$. When $m = n$ and
denoting $\mathbf{H}_{1}$ by $\mathbf{H}$, $(\mathbf{H}^{*}d\mathbf{H})$ is termed the
\emph{Haar measure} on $\mathfrak{U}^{\beta}(m)$.

The surface area or volume of the Stiefel manifold $\mathcal{V}^{\beta}_{m,n}$ is
\begin{equation}\label{vol}
    \Vol(\mathcal{V}^{\beta}_{m,n}) = \int_{\mathbf{H}_{1} \in
  \mathcal{V}^{\beta}_{m,n}} (\mathbf{H}^{*}_{1}d\mathbf{H}_{1}) =
  \frac{2^{m}\pi^{mn\beta/2}}{\Gamma^{\beta}_{m}[n\beta/2]},
\end{equation}
where $\Gamma^{\beta}_{m}[a]$ denotes the multivariate Gamma function for the space
$\mathfrak{S}_{m}^{\beta}$, and is defined by
\begin{eqnarray*}
  \Gamma_{m}^{\beta}[a] &=& \displaystyle\int_{\mathbf{A} \in \mathfrak{P}_{m}^{\beta}}
  \etr\{-\mathbf{A}\} |\mathbf{A}|^{a-(m-1)\beta/2 - 1}(d\mathbf{A}) \\
    &=& \pi^{m(m-1)\beta/4}\displaystyle\prod_{i=1}^{m} \Gamma[a-(i-1)\beta/2],
\end{eqnarray*}
where $\etr(\cdot) = \exp(\tr(\cdot))$, $|\cdot|$ denotes the determinant and $\re(a)
> (m-1)\beta/2$, see \citet{gr:87}.

Let $C_{\kappa}^{\beta}(\mathbf{B})$ be the Jack polynomials of $\mathbf{B} =
\mathbf{B}^{*}$, corresponding to the partition $\kappa=(k_{1},\ldots k_{m})$ of $k$,
$k_{1} \geq \cdots \geq k_{m} \geq 0$ with $\sum_{i=1}^{m}k_{i}=k$, see \citet{S:97}
and \citet{KE:06}. In addition,
$$
  {}_{p}F_{q}^{\beta}(a_{1}, \dots,a_{p};b_{1}, \dots,b_{q}; \mathbf{B}) = \sum_{k=0}^{\infty}\sum_{\kappa}
  \frac{[a_{1}]_{\kappa}^{\beta}, \dots, [a_{p}]_{\kappa}^{\beta}}{[b_{1}]_{\kappa}^{\beta}, \dots, [b_{q}]_{\kappa}}^{\beta}
  \frac{C_{\kappa}^{\beta}(\mathbf{B})}{k!},
$$
defines the hypergeometric function with one matrix argument on the space of hermitian
matrices, where $[a]^{\beta}_{\kappa}$ denotes the generalised Pochhammer symbol of weight
$\kappa$, defined as
$$
  [a]^{\beta}_{\kappa} = \prod_{i=1}^{m}(a-(i-1)\beta/2)_{k_{1}}
$$
where $ \Re(a) > (m-1)\beta/2 - k_{m}$ and $(a)_{i} = a(a+1)\cdots(a+i-1)$, see
\citet{gr:87}, \citet{KE:06} and \citet{dg:09}.

\begin{lemma}\label{lemj}
If $\mathbf{X} \in \mathfrak{L}^{\beta}_{m,n}$, then
\begin{equation}\label{eqj}
    \int_{\mathbf{H}_{1} \in \mathcal{V}_{m,n}^{\beta}} (\tr(\mathbf{XH}_{1}))^{2k} (d\mathbf{H}_{1}) =
  \sum_{\kappa}\frac{\left(
  \frac{1}{2}\right)_{k}}{[\beta n/2]^{\beta}_{\kappa}}C^{\beta}_{\kappa}(\mathbf{XX}^{*}),
\end{equation}
\end{lemma}
See \citet{dggj:10}.

Now, we use the complexification $\mathfrak{S}_{m}^{\beta, \mathfrak{C}} =
\mathfrak{S}_{m}^{\beta} + i \mathfrak{S}_{m}^{\beta}$ of $\mathfrak{S}_{m}^{\beta}$. That is,
$\mathfrak{S}_{m}^{\beta, \mathfrak{C}}$ consists of all matrices $\mathbf{X} \in
(\mathfrak{F^{\mathfrak{C}}})^{m \times m}$ of the form $\mathbf{Z} = \mathbf{X} +
i\mathbf{Y}$, with $\mathbf{X}, \mathbf{Y} \in \mathfrak{S}_{m}^{\beta}$. We refer to
$\mathbf{X} = \re(\mathbf{Z})$ and $\mathbf{Y} = \im(\mathbf{Z})$ as the \emph{real and
imaginary parts} of $\mathbf{Z}$, respectively. The \emph{generalised right half-plane}
$\mathbf{\Phi} = \mathfrak{P}_{m}^{\beta} + i \mathfrak{S}_{m}^{\beta}$ in
$\mathfrak{S}_{m}^{\beta,\mathfrak{C}}$ consists of all $\mathbf{Z} \in
\mathfrak{S}_{m}^{\beta,\mathfrak{C}}$ such that $\re(\mathbf{Z}) \in
\mathfrak{P}_{m}^{\beta}$, see \citet[p. 801]{gr:87}.

The next result generalises one given in  \citet{XuFang:89}, \citet{Teng89} and \citet{cdg:09}
for the real normed division algebras, see \citet{dg:09}:

\begin{lemma}\label{lem1}
Let $\mathbf{Z} \in \mathbf{\Phi}$ and $\mathbf{U} \in \mathfrak{S}_{m}^{\beta}$.
Then,
\begin{equation}\label{eqint}
    \int_{\mathbf{X}\in \mathfrak{P}_{m}^{\beta}}h(\tr \mathbf{XZ})|\mathbf{X}|^{a-(m-1)\beta/2-1}
    C_{\kappa}^{\beta}(\mathbf{XU})(d\mathbf{X}) =\frac{[a]_{\kappa}^{\beta}
    \Gamma_{m}^{\beta}[a]}{\Gamma[am+k]}|\mathbf{Z}|^{-a}C_{\kappa}^{\beta}(\mathbf{UZ}^{-1})\gamma,
\end{equation}
where
\begin{equation}\label{eqint2}
  \gamma=\int_{z \in \mathfrak{P}_{1}^{\beta}}h(z)z^{am+k-1}dw<\infty,
\end{equation}
for $\re(a)> (m-1)\beta/2 - k_{m}$, $\kappa=(k_{1},\ldots k_{m})$ a partition of $k$,
$k_{1} \geq \cdots \geq k_{m} \geq 0$ with $\sum_{i=1}^{m}k_{i}=k$.
\end{lemma}

\section{Affine shape distribution}\label{sec3}

Start with the  following   definition extended from \citet{gm:93}.

\begin{definition}\label{def1}
Two figures $\mathbf{X} \in \mathcal{L}^{\beta}_{K,N}$ and $\mathbf{X}_{1} \in
\mathcal{L}^{\beta}_{K,N}$ have the \textit{same configuration}, or \textit{affine shape}, if
$\mathbf{X}_{1}=\mathbf{XE}+\mathbf{1}_{N}e^{*}$, for some translation $e \in
\mathcal{L}^{\beta}_{1,N}$ and a $\mathbf{E} \in \mathcal{L}^{\beta}_{K,K}$.
\end{definition}

The configuration coordinates are constructed in the two steps summarised in the
expression
\begin{equation}\label{eqasd}
  \mathbf{LX}=\mathbf{Y}=\mathbf{UE}.
\end{equation}
The matrix $\mathbf{U} \in \mathcal{L}^{\beta}_{K,N-1}$ contains configuration coordinates of
$\mathbf{X}$. Let $\mathbf{Y}_{1} \in \mathcal{L}^{\beta}_{K,K}$ and $\mathbf{Y}_{2}  \in
\mathcal{L}^{\beta}_{K,q}$, with $q=N-K-1\geq 1$, such that $\mathbf{Y}=(\mathbf{Y}_{1}^{*}\mid
\mathbf{Y}_{2}^{*})^{*}$. Define also $\mathbf{U}=(\mathbf{I}\mid \mathbf{V}^{*})^{*}$, then
$\mathbf{V}=\mathbf{Y}_{2}\mathbf{Y}_{1}^{-1}$ and $\mathbf{E}=\mathbf{Y}_{1}$ where
$\mathbf{L}  \in \mathcal{L}^{\beta}_{N,N-1}$ is a Helmert sub-matrix.

Consider the following extension of \citet[Lemma 8]{cdg:09} for real normed division algebras.

\begin{lemma}\label{jacob}
Let $(\mathbf{F}^{1/2})^{2}=\mathbf{F} \in \mathfrak{P}_{K}^{\beta}$,  $\mathbf{H}
\in \mathfrak{U}^{\beta}(K)$, and let $\mathbf{E}=\mathbf{F}^{1/2} \mathbf{H}$ such
that $\mathbf{Y}=\mathbf{UF}^{1/2}\mathbf{H}$. Then
$$
  (d\mathbf{Y})=2^{-K}|\mathbf{F}|^{(q-1)\beta/2}(d\mathbf{V})(d\mathbf{F})(\mathbf{H}^{*}d\mathbf{H}).
$$
\end{lemma}
\begin{proof} Let $\mathbf{E}=\mathbf{\mathbf{F}}^{1/2}\mathbf{H}$, where $\mathbf{E} \in
\mathcal{L}_{K,K}^{\beta}$, $\mathbf{H} \in \mathfrak{U}^{\beta}(K)$ and $\mathbf{F}^{1/2} \in
\mathfrak{P}_{K}^{\beta}$. Therefore
$$
  \mathbf{E}^{*}\mathbf{E}=\mathbf{H}^{*}\mathbf{F}\mathbf{H}.
$$
Then from \citet[Lemma 2.2]{dggj:09}
$$
  (d(\mathbf{E}^{*}\mathbf{E}))=|\mathbf{H}^{*}\mathbf{H}|^{(K-1)\beta/2 +1}(d\mathbf{F})=
  |\mathbf{H}|^{(K-1)\beta + 2}(d\mathbf{F}) = (d\mathbf{F}).
$$
Now, by \citet[Lemma 2.15]{dggj:09}
\begin{eqnarray}
  (d\mathbf{E}) &=& 2^{-K}|\mathbf{E}^{*}\mathbf{E}|^{\beta/2-1}(d(\mathbf{E}^{*}\mathbf{E}))
        (\mathbf{H}^{*}d\mathbf{H}) \nonumber\\
    &=& 2^{-K}|\mathbf{H}^{*}\mathbf{F}\mathbf{H}|^{\beta/2-1}(d\mathbf{F})(\mathbf{H}^{*}d\mathbf{H}) \nonumber\\
    \label{eqjac1}
    &=& 2^{-K}|\mathbf{F}|^{\beta/2-1}(d\mathbf{F})(\mathbf{H}^{*}d\mathbf{H})
\end{eqnarray}
By other hand, observe that
$$
  \mathbf{Y}=\left(\begin{array}{c}
  \mathbf{I} \\
  \mathbf{V}
\end{array}\right)\mathbf{E}=\left(\begin{array}{c}
  \mathbf{E} \\
  \mathbf{V}\mathbf{E}
\end{array}\right).
$$
And by differentiating and computing the exterior product, see
\citet[Lemma 2.1]{dggj:09}, we obtain
$$
  (d\mathbf{Y})=|\mathbf{E}|^{\beta q}(d\mathbf{V})(d\mathbf{E}),
$$
but observing that $|\mathbf{E}|=|\mathbf{F}^{1/2}\mathbf{H}|=|\mathbf{F}|^{1/2}$, we
get
\begin{equation}\label{eqjac2}
 (d\mathbf{Y})=|\mathbf{F}|^{\beta q/2}(d\mathbf{V})(d\mathbf{E}).
\end{equation}
By replacing  (\ref{eqjac1}) into (\ref{eqjac2}) the desired result is obtained. \qed
\end{proof}

Now, recall that $\mathbf{X}\in \mathfrak{L}^{\beta}_{m,n}$ has a matrix multivariate
elliptically contoured distribution for real normed division algebras if its density, with
respect to the Lebesgue measure, is given by (see \citet{dggj:09}):
$$
  f_{\mathbf{X}}(\mathbf{X})=\frac{1}{|\mathbf{\Sigma}|^{\beta n/2}|\mathbf{\Theta}|^{\beta m/2}}
  h\left\{\tr\left[\mathbf{\Sigma}^{-1}(\mathbf{X}-\boldsymbol{\mu})^{*}\mathbf{\Theta}^{-1}
  (\mathbf{X}- \boldsymbol{\mu})\right]\right\},
$$
where  $\boldsymbol{\mu}\in \mathfrak{L}^{\beta}_{m,n}$, $ \mathbf{\Sigma}\in
\mathfrak{P}^{\beta}_{m}$,  $ \mathbf{\Theta}\in \mathfrak{P}^{\beta}_{m}$. The function $h:
\mathfrak{F} \rightarrow [0,\infty)$ is termed the generator function, and it is such that
$\int_{\mathfrak{P}^{\beta}_{1}} u^{\beta nm-1}h(u^2)du < \infty$.

Such a distribution is denoted by $\mathbf{X}\sim \mathcal{E}^{\beta}_{n\times
m}(\boldsymbol{\mu},\mathbf{\Sigma} \otimes \mathbf{\Theta}, h)$, for real case see
\citet{fz:90} and \citet{gv:93} and \citet{mdm:06} for complex case. Observe that this class of
matrix multivariate distributions includes gaussian, contaminated normal, Pearson type II and
VI, Kotz, Jensen-Logistic, power exponential, Bessel, among other distributions; whose
distributions have tails that are weighted more or less, and/or distributions with greater or
smaller degree of kurtosis than the gaussian distribution.

Now, we have the mathematical and statistics tools for establishing
the main density.
\begin{theorem}\label{teo1}
    Let $\mathbf{X}\sim \mathcal{E}_{N-1\times K}^{\beta}(\boldsymbol{\mu}_{\mathbf{X}},
    \mathbf{\Sigma}_{\mathbf{X}} \otimes \mathbf{\Theta},h)$. Then the affine shape density is given by
\begin{eqnarray}\label{ASD}
  \frac{\pi^{\beta K^{2}/2} \Gamma_{K}^{\beta}\left[\beta(N-1)/2\right]}{\Gamma_{K}^{\beta}\left[\beta K/2\right]
  |\mathbf{\Sigma}|^{\beta K/2} |\mathbf{U}^{*}\mathbf{\Sigma}^{-1}\mathbf{U}|^{\beta(N-1)/2}}
  \sum_{t=0}^{\infty}\frac{1}{t!\Gamma\left[K(N-1)/2+t\right]}\sum_{r=0}^{\infty}\frac{\left[
  \tr \mathbf{\Omega}\right]^{r}}{r!}\nonumber  \\
  \times \sum_{\tau}\frac{\left[\beta(N-1)/2\right]_{\tau}^{\beta}}{\left[\beta K/2\right]_{\tau}^{\beta}}
  C_{\tau}^{\beta}(\mathbf{U}^{*}\mathbf{\Omega}\mathbf{\Sigma}^{-1}\mathbf{U}(\mathbf{U}^{*}
  \mathbf{\Sigma}^{-1}\mathbf{U})^{-1})\gamma, \quad
\end{eqnarray}
where
\begin{equation}\label{eq: SCDnonisotropic}
  \gamma=\int_{z \in \mathfrak{P}_{1}^{\beta}}h^{(2t+r)}(z)z^{\beta K(N-1)2+t-1}dz <
  \infty,
\end{equation}
and $\mathbf{\Sigma} = \mathbf{L\Sigma L}^{*}$, $\boldsymbol{\mu} =
\boldsymbol{L\mu_{\mathbf{X}}}$ and
$\mathbf{\Omega}=\mathbf{\Sigma}^{-1}\boldsymbol{\mu}\mathbf{\Theta}\boldsymbol{\mu}^{*}$.
\end{theorem}
\begin{proof}
Define
$$
  \mathbf{LX\Theta}^{-1/2} = \mathbf{LZ} = \mathbf{Y} = \mathbf{UE},
$$
where $\left(\mathbf{\Theta}^{1/2}\right)^{2} = \mathbf{\Theta}$. Thus $\mathbf{Y}\sim
\mathcal{E}_{N-1\times K}^{\beta}(\boldsymbol{\mu\Theta}^{-1/2}, \mathbf{\Sigma} \otimes
\mathbf{I},h)$ where $\mathbf{\Sigma} = \mathbf{L\Sigma L}^{*}$, $\boldsymbol{\mu} =
\boldsymbol{L\mu_{\mathbf{X}}}$. Therefore the density of $\mathbf{Y}$ is given by
$$
  \frac{1}{|\mathbf{\Sigma}|^{\beta K/2}}\,h\{\tr[\mathbf{\Sigma}^{-1}(\mathbf{Y} -
  \boldsymbol{\mu\Theta}^{-1/2})(\mathbf{Y}-\boldsymbol{\mu\Theta}^{-1/2})^{*}]\}.
$$
Making the factorisation $\mathbf{Y} = \mathbf{UE} = \mathbf{U}\mathbf{F}^{1/2}\mathbf{H}$ and
by Lemma \ref{jacob}, then the joint density of $\mathbf{U}$, $\mathbf{F}$ and $\mathbf{H}$ is
\begin{eqnarray*}
  \frac{|\mathbf{F}|^{\beta(q+1)/2-1}}{2^{K}|\mathbf{\Sigma}|^{\beta K/2}}
  h\left[\tr\left(\mathbf{F}\mathbf{U}^{*}\mathbf{\Sigma}^{-1}\mathbf{U} + \mathbf{\Omega}\right)
  +\tr\left(-2\mathbf{\Theta}^{-1/2}\boldsymbol{\mu}^{*}\mathbf{\Sigma}^{-1}\mathbf{U}\mathbf{F}^{1/2}
  \mathbf{H}\right)\right]\hspace{1.5cm}\\
  \times \ (\mathbf{H}^{*}d\mathbf{H})(d\mathbf{F})(d\mathbf{V}),
\end{eqnarray*}
where
$\mathbf{\Omega}=\mathbf{\Sigma}^{-1}\boldsymbol{\mu}\mathbf{\Theta}^{-1}\boldsymbol{\mu}^{*}$.

Expanding $h$ in series of power, the joint density of $\mathbf{U}$, $\mathbf{F}$ and
$\mathbf{H}$ becomes:
\begin{eqnarray*}
  \frac{|\mathbf{F}|^{\beta(q+1)/2-1}}{2^{K}|\mathbf{\Sigma}|^{\beta K/2}}
  \sum_{t=0}^{\infty}\frac{1}{t!} h^{(t)}\left[\tr\left(\mathbf{F}\mathbf{U}^{*}
  \mathbf{\Sigma}^{-1}\mathbf{U }+  \mathbf{\Omega}\right)\right]
  \left[\tr\left(-2\mathbf{\Theta}^{-1/2}\boldsymbol{\mu}^{*}\mathbf{\Sigma}^{-1}\mathbf{U}
  \mathbf{F}^{1/2} \mathbf{H}\right)\right]^{t} \\
  (\mathbf{H}^{*}d\mathbf{H})(d\mathbf{F})(d\mathbf{V}).
\end{eqnarray*}
Now, using the Lemma \ref{lemj} for to integrate with respect to $\mathbf{H}$, and recalling
that $(1/2)_{t}4^{t}/(2t)! = 1/t!$ and $C_{\tau}^{\beta}(a\mathbf{B}) =
a^{t}C_{\tau}^{\beta}(\mathbf{B})$, the marginal joint density of $\mathbf{F}$ and $\mathbf{U}$
is
\begin{eqnarray*}
  \frac{\pi^{\beta K^{2}/2}|\mathbf{F}|^{\beta(q+1)/2-1}}{|\mathbf{\Sigma}|^{\beta K/2}
  \Gamma_{K}^{\beta}\left[\beta K/2\right]}
  \sum_{t=0}^{\infty}\frac{1}{t!}h^{(2t)}\left[\tr\left(\mathbf{F}\mathbf{U}^{*}
  \mathbf{\Sigma}^{-1}\mathbf{U} + \mathbf{\Omega}\right)\right]\hspace{4cm}\\
  \times \sum_{\tau}\frac{1}{\left[\beta K/2\right]_{\tau}^{\beta}}
  C_{\tau}^{\beta}\left(\mathbf{U}^{*}\mathbf{\Omega}
  \mathbf{\Sigma}^{-1}\mathbf{U}\mathbf{F}\right) (d\mathbf{F})(d\mathbf{V}),
\end{eqnarray*}
Assuming that $h^{2t}(\cdot)$ can be expanding in series of power, then the joint density of
$\mathbf{F}$ and $\mathbf{U}$ is

\begin{eqnarray*}
  \frac{\pi^{\beta K^{2}/2}|\mathbf{F}|^{\beta(q+1)/2-1}}{|\mathbf{\Sigma}|^{\beta K/2}
  \Gamma_{K}^{\beta}\left[\beta K/2\right]}
  \sum_{t=0}^{\infty}\frac{1}{t!}\sum_{r=0}^{\infty}\frac{1}{r!}
  h^{(2t+r)}\left[\tr \left(\mathbf{F}\mathbf{U}^{*} \mathbf{\Sigma}^{-1} \mathbf{U} \right) \right]
  \left[\tr\left(\boldsymbol{\mu}^{*}\mathbf{\Sigma}^{-1}\boldsymbol{\mu}\right)\right]^{r}\hspace{.5cm}\\
  \times \sum_{\tau}\frac{1}{\left[\beta K/2\right]_{\tau}^{\beta}}
  C_{\tau}^{\beta}\left(\mathbf{U}^{*}\mathbf{\Omega}
  \mathbf{\Sigma}^{-1}\mathbf{U}\mathbf{F}\right) (d\mathbf{F})(d\mathbf{V}).
\end{eqnarray*}
Hence, the marginal density of $\mathbf{U}$ is
\begin{eqnarray}\label{eqad}
  \frac{\pi^{\beta K^{2}/2}}{|\mathbf{\Sigma}|^{\beta K/2}
  \Gamma_{K}^{\beta}\left[\beta K/2\right]}
  \sum_{t=0}^{\infty}\frac{1}{t!}\sum_{r=0}^{\infty}\frac{1}{r!}\left[
  \tr\left(\boldsymbol{\mu}^{*}\mathbf{\Sigma}^{-1}\boldsymbol{\mu}\right)\right]^{r}
  \sum_{\tau}\frac{1}{\left[\beta K/2\right]_{\tau}^{\beta}}\nonumber \hspace{2cm}\\
  \times \int_{\mathbf{F}>\mathbf{0}}h^{(2t+r)}
  \left(\tr\left(\mathbf{F}\mathbf{U}^{*}\mathbf{\Sigma}^{-1}\mathbf{U}\right)\right)
  |\mathbf{F}|^{\beta(q+1)/2-1} C_{\tau}^{\beta}\left(\mathbf{U}^{*}\mathbf{\Omega}\mathbf{\Sigma}^{-1}\mathbf{U}
  \mathbf{F}\right)(d\mathbf{F}).
\end{eqnarray}
From (\ref{eqint}), the integral in (\ref{eqad}) is evaluated as
\begin{eqnarray*}
\int_{\mathbf{F}>\mathbf{0}}h^{(2t+r)}
  \left[\tr\left(\mathbf{F}\mathbf{U}^{*}\mathbf{\Sigma}^{-1}\mathbf{U}\right)\right]|\mathbf{F}|^{\beta(q+1)/2-1}
  C_{\tau}^{\beta}\left(\mathbf{U}^{*}\mathbf{\Omega}\mathbf{\Sigma}^{-1}\mathbf{U}\mathbf{F}\right)
  (d\mathbf{F})\hspace{2cm}\\
  =\frac{\left[\beta(N-1)/2\right]_{\tau}\Gamma_{K}^{\beta}\left[\beta(N-1)/2\right]}
  {\Gamma\left[K(N-1)/2+t\right]}
  \frac{C_{\tau}^{\beta}(\mathbf{U}^{*}\mathbf{\Omega}\mathbf{\Sigma}^{-1}\mathbf{U}(\mathbf{U}^{*}
  \mathbf{\Sigma}^{-1}\mathbf{U})^{-1})}{|\mathbf{U}^{*}\mathbf{\Sigma}^{-1}
  \mathbf{U}|^{\beta(N-1)/2}}\gamma,
\end{eqnarray*}
where
$$
  \gamma=\int_{z \in \mathfrak{P}_{1}^{\beta}}h^{(2t+r)}(z)z^{\beta K(N-1)/2+t-1}dz<\infty,
$$
and the required result follows. \qed
\end{proof}

Note that the general density is indexed by a simple univariate integral involving the general
derivative of generator function. These kind of densities appear rare in matrix-variate
distributions (see \citet{cdg:09}). However, they demand the computation of derivatives of any
order, which is not a trivial fact; general formulae for the classical elliptical models (Kotz,
Pearson, Bessel, Jensen-Logistic) are available too   in the above mentioned reference.

Finally, the central and isotropic affine shape densities are obtained.

\begin{corollary}\label{coro1}
If $\mathbf{X}\sim \mathcal{E}_{N-1\times K}^{\beta}(\mathbf{0}, \mathbf{\Sigma}_{\mathbf{X}}
\otimes \mathbf{\Theta},h)$, then the central affine shape density is invariant under the
elliptical distributions, moreover, its density is
$$
  \frac{\Gamma_{K}^{\beta}\left[\beta(N-1)/2\right]}{\pi^{\beta Kq/2}
  \Gamma_{K}^{\beta}\left[\beta K/2\right]|\mathbf{\Sigma}|^{\beta K/2}
  }|\mathbf{U}^{*}\mathbf{\Sigma}^{-1}\mathbf{U}|^{-\beta (N-1)/2}.
$$
\end{corollary}
\begin{proof} The proof  follows by taking $t= r =0$ in (\ref{ASD}) and noting that $h^{(2t+r)}(z) = h^{(0)}(z) \equiv
h(z)$, hence
$$
  \frac{\pi^{\beta K^{2}/2}\Gamma_{K}^{\beta}\left[\beta(N-1)/2\right]}{|\mathbf{\Sigma}|^{\beta K/2}
  \Gamma_{K}^{\beta}\left[\beta K/2\right]}\frac{|\mathbf{U}^{*}\mathbf{\Sigma}^{-1}\mathbf{U}|^{-\beta(N-1)/2}}
  {\Gamma\left[\beta K(N-1)/2\right]}\gamma.
$$
Now, using \citet[p. 59]{fz:90},
$$
  \gamma =  \int_{0}^{\infty}h(z)z^{\beta K(N-1)/2-1}dz =
  \frac{\Gamma\left[\beta K(N-1)/2\right]}{\pi^{\beta K(N-1)/2}},
$$
the desired result is obtained. \qed
\end{proof}

\begin{corollary}\label{ASDiso}
If\, $\mathbf{Y}\sim \mathcal{E}_{N-1\times K}^{\beta}(\boldsymbol{\mu}_{\mathbf{X}},
\sigma^{2}\mathbf{I}_{N-1}\otimes \mathbf{\Theta},h)$,  then the isotropic noncentral affine
shape density is given by
\begin{eqnarray*}
  \frac{\pi^{\beta K^{2}/2}\Gamma_{K}^{\beta}\left[\beta(N-1)/2\right)]}{\Gamma_{K}^{\beta}
  \left[\beta K/2\right]|\mathbf{I}_{K} + \mathbf{V}^{*}\mathbf{V}|^{\beta(N-1)/2}}
  \sum_{t=0}^{\infty}\frac{1}{t!\Gamma\left[K(N-1)/2+t\right]}\sum_{r=0}^{\infty}\frac{1}{r!}\left[
  \tr\left(\mathbf{\Omega}_{1}\right)\right]^{r}\hspace{1cm}\\
  \times \sum_{\tau}\frac{\left[\beta(N-1)/2\right]_{\tau}^{\beta}}{\left[\beta K/2\right]_{\tau}^{\beta}}
  C_{\tau}^{\beta}\left(\mathbf{U}^{*}\boldsymbol{\Omega}_{1}\mathbf{U}(\mathbf{I}_{K} +
  \mathbf{V}^{*}\mathbf{V})^{-1}\right)\gamma,
\end{eqnarray*}
where
$$
  \gamma=\int_{z \in \mathfrak{P}_{1}^{\beta}}h^{(2t+r)}(z)z^{\beta K(N-1)/2 + t -1}dz<\infty,
$$
and $\mathbf{\Omega}_{1} = \sigma^{-2} \boldsymbol{\mu} \mathbf{\Theta}^{-1}
\boldsymbol{\mu}^{*}$.
\end{corollary}
\begin{proof}
The result follows easily, just recall  that
$\mathbf{U}=(\mathbf{I}_{K}\mid \mathbf{V}^{*})^{*}$ and take
$\mathbf{\Sigma} = \sigma^{2}\mathbf{I}_{N-1}$ in (\ref{ASD}). Thus
$$
  |\mathbf{\Sigma}|^{\beta K/2} |\mathbf{U}^{*}\mathbf{\Sigma}^{-1}\mathbf{U}|^{\beta(N-1)/2} =
  |\mathbf{I}_{K} + \mathbf{V}^{*}\mathbf{V}|^{\beta(N-1)/2}. \mbox{\qed}
$$
\end{proof}

\section{The Gaussian affine shape distribution}\label{sec4}

In this section we study the Gaussian affine shape distribution, as  corollary of the preceding
results.
\begin{corollary}\label{cor41}
    Let $\mathbf{X}\sim \mathcal{E}_{N-1\times K}^{\beta}(\boldsymbol{\mu}_{\mathbf{X}},
    \mathbf{\Sigma}_{\mathbf{X}} \otimes \mathbf{\Theta},h)$. Then the affine shape density is given by
\begin{eqnarray}\label{ASDG}
  \frac{ \Gamma_{K}^{\beta}\left[\beta(N-1)/2\right] \etr\{-\beta (\mathbf{\Omega} - \mathbf{U}^{*}
  \mathbf{\Omega}\mathbf{\Sigma}^{-1}\mathbf{U}(\mathbf{U}^{*}
  \mathbf{\Sigma}^{-1}\mathbf{U})^{-1})/2\}}{\pi^{\beta Kq/2}\Gamma_{K}^{\beta}\left[\beta K/2\right]
  |\mathbf{\Sigma}|^{\beta K/2} |\mathbf{U}^{*}\mathbf{\Sigma}^{-1}\mathbf{U}|^{\beta(N-1)/2}}
  \nonumber  \\
  {}_{1}F_{1}^{\beta}(-\beta q/2; \beta K/2; -\beta \mathbf{U}^{*}\mathbf{\Omega}
  \mathbf{\Sigma}^{-1}\mathbf{U}(\mathbf{U}^{*} \mathbf{\Sigma}^{-1}\mathbf{U})^{-1}/2).
\end{eqnarray}
\end{corollary}
\begin{proof}
From \citet{dggj:09}, the gaussian case turns,
$$
  h(v)= \frac{1}{\displaystyle \left(\frac{2\pi}{\beta}\right)^{\beta K(N-1)/2}} \exp\{-\beta v/2\}
$$
and
$$
  h^{2t+r}(v)= \frac{\displaystyle \left(-\frac{\beta}{2}\right)^{2t+r}}
  {\displaystyle \left(\frac{2\pi}{\beta}\right)^{\beta K(N-1)/2}} \exp\{-\beta v/2\}.
$$
Therefore
\begin{eqnarray*}
  \gamma &=& \int_{z \in \mathfrak{P}_{1}^{\beta}}h^{(2t+r)}(z)z^{\beta K(N-1)/2 + t -1}dz \\
    &=& \int_{z \in \mathfrak{P}_{1}^{\beta}}\frac{\displaystyle \left(-\frac{\beta}{2}\right)^{2t+r}}
  {\displaystyle \left(\frac{2\pi}{\beta}\right)^{\beta K(N-1)/2}} \exp\{-\beta z/2\} dz \\
    &=& \displaystyle\frac{\Gamma[\beta K(N-1)/2 + t]}{\pi^{\beta K(N - 1)/2}}
    \left(-\frac{\beta}{2}\right)^{r}\left(\frac{\beta}{2}\right)^{t}.
\end{eqnarray*}
Also, note that
$$
  [\tr \mathbf{\Omega}]^{r}  \left(-\frac{\beta}{2}\right)^{r} = [\tr \mathbf{-\beta \Omega}/2]^{r}
$$
and
$$
  \left(\frac{\beta}{2}\right)^{t}C_{\tau}^{\beta}(\mathbf{U}^{*}\mathbf{\Omega}
  \mathbf{\Sigma}^{-1}\mathbf{U}(\mathbf{U}^{*}\mathbf{\Sigma}^{-1}\mathbf{U})^{-1}) =
  C_{\tau}^{\beta}(\beta\mathbf{U}^{*}\mathbf{\Omega}\mathbf{\Sigma}^{-1}\mathbf{U}(\mathbf{U}^{*}
  \mathbf{\Sigma}^{-1}\mathbf{U})^{-1}/2).
$$
The final expression for (\ref{ASDG}) is obtained by applying the Kummer relations, see
\citet{dg:09}. \qed
\end{proof}

The isotropic case of this distribution is given by

\begin{corollary}\label{cor411}
    Let $\mathbf{X}\sim \mathcal{E}_{N-1\times
K}^{\beta}(\boldsymbol{\mu}_{\mathbf{X}},
    \mathbf{\Sigma}_{\mathbf{X}} \otimes \mathbf{\Theta},h)$. Then the affine shape density is given by
\begin{eqnarray}\label{ASDGI}
  \frac{ \Gamma_{K}^{\beta}\left[\beta(N-1)/2\right] \etr\{-\beta (\mathbf{\Omega} - \mathbf{U}^{*}
  \mathbf{\Omega}\mathbf{\Sigma}^{-1}\mathbf{U}(\mathbf{U}^{*}
  \mathbf{\Sigma}^{-1}\mathbf{U})^{-1})/2\}}{\pi^{\beta Kq/2}\Gamma_{K}^{\beta}\left[\beta K/2\right]
  |\mathbf{I}_{K} + \mathbf{V}^{*}\mathbf{V}|^{\beta(N-1)/2}}
  \nonumber  \\
  {}_{1}F_{1}^{\beta}(-\beta q/2; \beta K/2; -\beta \mathbf{U}^{*}\mathbf{\Omega}
  \mathbf{\Sigma}^{-1}\mathbf{U}(\mathbf{U}^{*} \mathbf{\Sigma}^{-1}\mathbf{U})^{-1}/2)\\
    =\frac{ \Gamma_{K}^{\beta}\left[\beta(N-1)/2\right] \etr\{-\beta (\sigma^{-2}\boldsymbol{\mu}\boldsymbol{\mu}^{*} - \sigma^{-2}\mathbf{U}^{*}
  \boldsymbol{\mu}\boldsymbol{\mu}^{*}\mathbf{U}(\mathbf{U}^{*}
  \mathbf{U})^{-1})/2\}}{\pi^{\beta Kq/2}\Gamma_{K}^{\beta}\left[\beta K/2\right]
  |\mathbf{I}_{K} + \mathbf{V}^{*}\mathbf{V}|^{\beta(N-1)/2}}
  \nonumber  \\
  {}_{1}F_{1}^{\beta}(-\beta q/2; \beta K/2; -\beta\sigma^{-2} \mathbf{U}^{*}\boldsymbol{\mu}\boldsymbol{\mu}^{*}
  \mathbf{U}(\mathbf{U}^{*} \mathbf{U})^{-1}/2).
\end{eqnarray}
\end{corollary}

It follows straightforwardly, recall that $\mathbf{U}=(\mathbf{I}_{K}\mid \mathbf{V}^{*})^{*}$
and replace $\mathbf{\Sigma} = \sigma^{2}\mathbf{I}_{N-1}$ in (\ref{ASDGI}). Thus
$$
  |\mathbf{\Sigma}|^{\beta K/2} |\mathbf{U}^{*}\mathbf{\Sigma}^{-1}\mathbf{U}|^{\beta(N-1)/2} =
  |\mathbf{I}_{K} + \mathbf{V}^{*}\mathbf{V}|^{\beta(N-1)/2}. \mbox{\qed}
$$

\section{Example}\label{sec5}

In this last section we apply the complex normal shape model in a classical data, the brain
magnetic resonance scans of normal and schizophrenic patients (see \citet{bo:66}, \citet{DM98},
among many others). The random samples consist of  14  scans of each group, which correspond to
a  near midsagittal two dimensional slices of MR. On each image, an expert locates the
following 13 landmarks and registers their coordinates. The selected points correspond to
(\citet{bo:66}, \citet{DM98}): 1. splenium, 2. genu, 3. top of corpus callosum, 4. top of head,
5. tentorium of cerebellum at dura, 6. top of cerebellum, 7. tip of fourth ventricle, 8. bottom
of cerebellum, 9. top of pons, 10. bottom of pons, 11. optic chiasm, 12. frontal pole, 13.
superior colliculus (see \citet[figure 9, p.12]{DM98}). The aim of the experiment is to test
the equality in mean shape of the two population after removing the non geometrical information
of the scans.

The preceding works have analysed the shape under Euclidian transformation by filtering out
translation, scaling and rotation of scans, and concluding that both mean shapes are
statistically different. However, the shape theory via similarity transformation works well
when the objects are rigid and develop a``constant radial growth" in some sense, but in the
case of the brain, which is a soft organ, it is prone to deformations, so in order to match the
scans we need to filter out the shear, instead of the rotation. Then the affine shape or
configuration analysis is more appropriate than the usual Euclidian shape (this discrepancy can
be solved statistically by using the so termed modified BIC criteria, see \citet{YY07} and the
references therein).

Now, if the original landmark $13\times 1$ vector $\mathbf{X}$ follows an isotropic complex
normal model, then from corollary \ref{cor411} we have that the corresponding affine shape
density is given by
\begin{eqnarray}\label{ASDCG}
  \frac{ \Gamma_{K}^{\beta}\left[\beta(N-1)/2\right] \etr\{-\beta\sigma^{-2}\boldsymbol{\mu}\boldsymbol{\mu}^{*}/2 +\beta \sigma^{-2}\mathbf{U}^{*}
  \boldsymbol{\mu}\boldsymbol{\mu}^{*}\mathbf{U}(\mathbf{U}^{*}
  \mathbf{U})^{-1}/2\}}{\pi^{\beta Kq/2}\Gamma_{K}^{\beta}\left[\beta K/2\right]
  |\mathbf{U}^{*}\mathbf{U}|^{\beta(N-1)/2}}
  \nonumber  \\
  {}_{1}F_{1}^{\beta}(-\beta q/2; \beta K/2; -\beta\sigma^{-2} \mathbf{U}^{*}\boldsymbol{\mu}\boldsymbol{\mu}^{*}
  \mathbf{U}(\mathbf{U}^{*} \mathbf{U})^{-1}/2).
\end{eqnarray}

where $\beta=2$, $K=1$, $N=3$, $q=11$. Then it is of interest  the estimation of the scale
parameter $\sigma^{2}$ and the location parameter
\begin{eqnarray*}
\boldsymbol{\mu}=(\mu_{1,1}+i\mu_{1,2},\mu_{2,1}+i\mu_{2,2},\mu_{3,1}+i\mu_{3,2},\mu_{4,1}+i\mu_{4,2},\mu_{5,1}+i\mu_{5,2},
\mu_{6,1}+i\mu_{6,2},\\\mu_{7,1}+i\mu_{7,2},
\mu_{8,1}+i\mu_{8,2},\mu_{9,1}+i\mu_{9,2},\mu_{10,1}+i\mu_{10,2},\mu_{11,1}+i\mu_{11,2},\mu_{12,1}+i\mu_{12,2})'
\end{eqnarray*}

Note that the above affine density is a polynomial of degree eleven, then the inference
procedure is notoriously simple.

Let $\mathfrak{L}(\widetilde{\boldsymbol{\mu}},\widetilde{\sigma}^{2},h)$ be the log likelihood
function of a given group-model. The maximisation of the likelihood function
$\mathfrak{L}(\widetilde{\boldsymbol{\mu}},\widetilde{\sigma}^{2},h)$, is obtained in this
paper by using the \emph{Nelder-Mead Simplex Method}, which is an unconstrained multivariable
function using a derivative-free method; specifically, we apply  the routine
\textbf{fminsearch} implemented by the sofware MatLab.

The shape densities are polynomials of scalar zonal polynomials, i.e. hypergeometric series
which terminates and this can be computed easily by the algorithms of  \citet{KE:06}.

At this point the log likelihood can be computed, then we use fminsearch for the MLE's. The
initial value for the algorithm is the sample mean of the normal  variables
$\mathbf{Y}\sim\mathcal{N}_{12 \times 1}(\boldsymbol{\mu}, \sigma^{2}\mathbf{I}_{12}) $ and the
median of the variances.

The maximum likelihood estimators for location parameters associated with the normal and
schizophrenic groups under the complex Gaussian model, are the following:

For the normal group:
\begin{footnotesize}
\begin{eqnarray*}
\boldsymbol{\widetilde{\mu}}=(-1.5312+i0.6468,-0.0084+i0.4637,1.4306+i0.9779,1.2459-i2.1091,0.3962-i0.7200,\\
-0.0475-i1.3257,-0.2321-i1.9199,
-0.6937-i0.2605,-0.6556-i1.0118,-1.3664+i0.1697,\\-2.4834+i1.4590,0.2667+i0.0804)',
\end{eqnarray*}
\end{footnotesize} and $\widetilde{\sigma}^{2}=0.0159$  and for the
schizophrenic group:
\begin{footnotesize}
\begin{eqnarray*}
\boldsymbol{\widetilde{\mu}}=(-1.6898+i0.0545,0.0179+i0.5385,1.1107+i1.3683,1.7652-i1.7264,0.5794-i0.6484,\\
0.3610-i1.300,0.3015-i1.9547,
-0.6052-i0.4795,-0.3432-i1.2216,-1.3368-i0.2235,\\-2.8060+i0.6043,0.1087-i0.0434)'
\end{eqnarray*}
\end{footnotesize}
and $\widetilde{\sigma}^{2}=0.0155$.

Finally, we can test equality in affine shape between the two independent populations. In this
experiment we have: two independent samples of 14 patients and 24 population shape parameters
to estimate for each group. Namely, if $L(\boldsymbol{\mu}_{n},\boldsymbol{\mu}_{s})$ is the
likelihood, where $\boldsymbol{\mu}_{n}$, $\boldsymbol{\mu}_{s}$,  represent the mean shape
parameters of the normal and schizophrenic group, respectively, then we want to test:
$H_{0}:\boldsymbol{\mu}_{n}=\boldsymbol{\mu}_{s}$ vs
$H_{1}:\boldsymbol{\mu}_{n}\neq\boldsymbol{\mu}_{s}$. Then $-2\log \Lambda=2\sup _{H_{1}}\log
L(\boldsymbol{\mu}_{n},\boldsymbol{\mu}_{s}) -2\sup _{H_{0}}\log L(\boldsymbol{\mu}_{n},
\boldsymbol{\mu}_{s})$, and according to Wilk's theorem $-2\log \Lambda\sim\chi^{2}_{24}$ under
$H_{0}.$

In this case  we  obtained  that:
$$
  -2\log \Lambda=2(858.09)-2(834.60)=46.98,
$$
Since the p-value for the test is
$$
  P(\chi^{2}_{24}\geq 46.98)=0.0034
$$
we have important evidence that the normal and schizophrenic brain MR are different in affine
shape. This conclusion is ratified by \citet{DM98}, for example, but using a ``rigid" Euclidian
match in $\Re^{2}$, in this case they obtained a p-value of $0.005$.

Some new models (Kotz, Pearson VII, Bessel, Jensen-Logistic) can be studied and contrast  them
with the Gaussian one via BIC modified criteria, however new extensions of the so termed
generalised Kummer relations with Jack polynomials are required (a generalisation to real
normed division algebras of \citet{h:55}); this shall constitute part of a future work.

Finally, observe that the real dimension of real normed division algebras can be expressed as
potentia of 2, $\beta = 2^{n}$ for $n = 0,1,2,3$. On the other hand, as observed by
\citet{k:84}, the results obtained in this work can be extended to hypercomplex cases; that is,
for complex, bicomplex, biquaternion and bioctonion (or sedenionic) algebras, which of course
are not division algebras (except the complex algebra), but are Jordan algebras, and all their
isomorphic algebras. Note, too, that hypercomplex algebras are obtained by replacing the real
numbers with complex numbers in the construction of real normed division algebras. Thus, the
results for hypercomplex algebras are obtained by simply replacing $\beta$ with $2\beta$ in our
results (we reiterate, as reported by \citet{k:84}). Alternatively, following \citet{k:84}, we
can conclude that, our results are true for `$2^{n}$-ions', $n = 0,1,2,3,4,5$, emphasising that
only for $n=0,1,2,3$ are the result algebras in fact real normed division algebras.

\section*{Acknowledgments}
This research work was supported  by University of Medellin (Medellin, Colombia) and
Universidad Aut\'onoma Agraria Antonio Narro (M\'{e}xico),  joint grant No. 469, SUMMA group.
Also, the first author was partially supported  by IDI-Spain, Grants No. \ FQM2006-2271 and
MTM2008-05785 and the paper was written during J. A. D\'{\i}az- Garc\'{\i}a's stay as a
visiting professor at the Department of Statistics and O. R. of the University of Granada,
Spain.

\end{document}